# Resolution and simplification of Dombi-fuzzy relational equations and latticized optimization programming on Dombi FREs

**Amin Ghodousian [1], Sara Zal [2]**

[1]- First School of Engineering Science, College of Engineering
University of Tehran, P.O.Box 11365-4563, Tehran, Iran
Email: a.ghodousian@ut.ac.ir

[2]- Second School of Engineering Science, College of Engineering
University of Tehran, Tehran, Iran
Email: sara.zal@ut.ac.ir

**Abstract**

In this paper, we introduce a type of latticized optimization problem whose objective function is the maximum component function and the feasible region is defined as a system of fuzzy relational equalities (FRE) defined by the Dombi t-norm. Dombi family of t-norms includes a parametric family of continuous strict t-norms, whose members are increasing functions of the parameter. This family of t-norms covers the whole spectrum of t-norms when the parameter is changed from zero to infinity. Since the feasible solutions set of FREs is non-convex and the finding of all minimal solutions is an NP-hard problem, designing an efficient solution procedure for solving such problems is not a trivial job.

Some necessary and sufficient conditions are derived to determine the feasibility of the problem. The feasible solution set is characterized in terms of a finite number of closed convex cells. An algorithm is presented for solving this nonlinear problem. It is proved that the algorithm can find the exact optimal solution and an example is presented to illustrate the proposed algorithm.

1. **INTRODUCTION** (with two 9 pt lines space from the keywords)

In this paper, we study the following linear problem in which the constraints are formed as fuzzy relational equations defined by Dombi t-norm:

$$\min \quad z(x) = \max\{x_1, x_2, \ldots, x_n\}$$
$$A \varphi x = b$$
$$x \in [0,1]^n \tag{1}$$

Where $I = \{1, \ldots, m\}$, $J = \{1, \ldots, n\}$, $A = (a_{ij})_{m \times n}$ is a fuzzy matrix such that $0 \leq a_{ij} \leq 1$ ($\forall i \in I$ and $\forall j \in J$), $b = (b_i)_{m \times 1}$ is a fuzzy vector such that $0 \leq b_i \leq 1$ ($\forall i \in I$), and "$\varphi$" is the Dombi t-norm defined as follows:

$$\varphi(x,y) = \begin{cases} 0 & x = 0 \text{ or } y = 0 \\ \left(1 + \left[\left(\frac{1-x}{x}\right)^\lambda + \left(\frac{1-y}{y}\right)^\lambda\right]^{\frac{1}{\lambda}}\right)^{-1} & \text{otherwise} \end{cases} \tag{2}$$

where $\lambda > 0$. If $a_i$ is the $i$'th row of matrix $A$, then problem (1) can be expressed as follows:

$$\begin{aligned} min \quad & z(x) = \max\{x_1, x_2, \ldots, x_n\} \\ & \varphi(a_i. x) = b_i \quad . i \in I \\ & x \in [0.1]^n \end{aligned} \tag{3}$$

where the constraints mean

$$\varphi(a_i. x) = \max_{j \in J}\{\varphi(a_{ij}. x_j)\} = b_i \ (\forall i \in I) \text{ and}$$

$$\varphi(a_{ij}, x_j) = \begin{cases} 0 & a_{ij} = 0 \text{ or } x_j = 0 \\ \left(1 + \left[\left(\frac{1-a_{ij}}{a_{ij}}\right)^\lambda + \left(\frac{1-x_j}{x_j}\right)^\lambda\right]^{\frac{1}{\lambda}}\right)^{-1} & \text{otherwise} \end{cases} \tag{4}$$

As mentioned, the family $\varphi$ is increasing in $\lambda$. On the other hand, Dombi t-norm $\varphi$ converges to the basic fuzzy intersection $\min\{x, y\}$ as $\lambda$ goes to infinity and converges to Drastic product t-norm as $\lambda$ approaches zero. Therefore, Dombi t-norm covers the whole spectrum of t-norms [5]. In [3], the Dombi operations of single-valued neutrosophic numbers (SVNNs) were presented based on the operations of the Dombi t-norm and t-conorm. The authors proposed the single-valued neutrosophic Dombi weighted arithmetic average operator and the single-valued neutrosophic Dombi weighted geometric average operator to deal with the aggregation of SVNNs. In [1], a fuzzy morphological approach was presented to detect the edges of real time images in order to preserve their features, where Dombi t-norm and t-conorm was used for computing morphological dilation and erosion. In [4], the authors studied the connection with Dombi aggregative operators, uninorms, strict t-norms and t-conorms. They presented a new representation theorem of strong negations that explicitly contains the neutral value.

The theory of fuzzy relational equations (FRE) was firstly proposed by Sanchez and applied in problems of the medical diagnosis [17]. Nowadays, it is well known that many issues associated with a body knowledge can be treated as FRE problems [16]. The solution set of FRE is often a non-convex set that is completely determined by one maximum solution and a finite number of minimal solutions [8]. The other bottleneck is concerned with detecting the minimal solutions that is an NP-hard problem [9,11,12,14]. The problem of optimization subject to FRE and FRI is one of the most interesting and on-going research topic among the problems related to FRE and FRI theory [2,6,8–12,15,18]. Recently, many interesting generalizations of the linear programming subject to a system of fuzzy relations have been introduced and developed based on composite operations used in FRE, fuzzy relations used in the definition of the constraints, some developments on the objective function of the problems and other ideas [7,11-13,15].

Yang et al. [19] studied the single-variable term semi-latticized geometric programming subject to max-product fuzzy relation equations. The proposed problem was devised from the peer-to-peer network system and the target was to minimize the biggest dissatisfaction degrees of the terminals in such system. Yang et al. [20] introduced another version of the latticized programming problem subject to max-prod fuzzy relation inequalities with application in the optimization management model of wireless communication emission base stations.

The latticized problem was defined by minimizing objective function z(x)=max{$x_1, x_2, \ldots, x_n$} subject to feasible region $X(A, b) = \{x \in [0,1]^n : A \circ x \geq b\}$ where "$\circ$" denotes fuzzy max-product composition.

The rest of the paper is organized as follows. In Section 2, a necessary and sufficient condition is derived to determine the feasibility of max-Dombi FRE. In Section 3, the feasible solution set of problem (1) is characterized. It is shown that the feasible region can be expressed as the union of a finite closed convex cells. Section 4 describes the optimal solution of Problem (1). An algorithm is proposed to find the optimal solution and finally, Section 5 provides a numerical example to illustrate the algorithm.

## 2. BASIC PROPERTIES OF MAX-DOMBI FRE REVIEW STAGE

Let $S$ denote the feasible solutions set of problem (1), that is, $S = \{x \in [0,1]^n : \max_{j=1}^{n}\{\varphi(a_{ij}, x_j)\} = b_i, \forall i \in I\}$. Also, for each $i \in I$, define $J_i = \{j \in J : a_{ij} \geq b_i\}$. Also, for each $i \in I$ and $j \in J_i$, we define

$$V(b_i, a_{ij}) = \left(1 + \left[\left(\frac{1-b_i}{b_i}\right)^\lambda - \left(\frac{1-a_{ij}}{a_{ij}}\right)^\lambda\right]^{\frac{1}{\lambda}}\right)^{-1}$$

(5)

According to [9], when $S \neq \emptyset$, it can be completely determined by one maximum solution and a finite number of minimum solutions. The maximum solution can be obtained by $\bar{X} = \min_{i \in I}\{\hat{x}_i\}$ where $\hat{x}_i = [(\hat{x}_i)_1, \dots, (\hat{x}_i)_n]$ ($\forall i \in I$) is defined as follows

$$(\hat{x}_i)_j = \begin{cases} V(b_i, a_{ij}) & j \in J_i, \ b_i \neq 0 \\ 0 & j \in J_i, \ a_{ij} > b_i = 0 \\ 1 & \text{otherwise} \end{cases}, \forall j \in J$$

(6)

Moreover, if we denote the set of all minimum solutions by $\underline{S}$, then

$$S = \bigcup_{\underline{X} \in \underline{S}}\{x \in [0,1]^n : \underline{X} \leq x \leq \bar{X}\}$$

(7)

**Corollary 1.** $S \neq \emptyset$ if and only if $\bar{X} \in S$.

**Definition 1.** Let $i \in I$. For each $j \in J_i$, we define $\check{x}_i(j) = [\check{x}_i(j)_1, \dots, \check{x}_i(j)_n]$ such that

$$\check{x}_i(j)_k = \begin{cases} V(b_i, a_{ik}) & b_i \neq 0 \text{ and } k = j \\ 0 & \text{otherwise} \end{cases}, \forall k \in J$$

(8)

**Definition 2.** Let $e: I \to \bigcup_{i \in I} J_i$ so that $e(i) = j \in J_i, \forall i \in I$, and let $E$ be the set of all vectors $e$. For the sake of convenience, we represent each $e \in E$ as an $m-$ dimensional vector $e = [j_1, \dots, j_m]$ in which $J_k = e(k)$.

**Definition 3.** Let $e \in E$. We define $\underline{X}(e) = [\underline{X}(e)_1, \dots, \underline{X}(e)_n]$ where $\underline{X}(e)_j = \max_{i \in I}\{\check{x}_i(e(i))_j\}, \forall j \in J$.

Based on the definitions 2 and 3, and according to [9], we have

$$\underline{S} \subseteq \{\underline{X}(e) : e \in E\}$$

(9)

Moreover, we have the following fundamental theorem [9]:

**Theorem 1.** Suppose that matrix $\bar{A} = (\bar{a}_{ij})_{m \times n}$ is resulted from $A = (a_{ij})_{m \times n}$ by the following steps:

(I) If $j_0 \notin J_i$ for some $i \in I$ and $j_0 \in J$, then we set $\bar{a}_{ij_0} = 0$.

(II) If $j_0 \in J_i$, $b_i \neq 0$ and there exists some $i' \in I$ ($i' \neq i$) such that $j_0 \in J_{i'}$, $b_{i'} \neq 0$ and $V(b_{i'}, a_{i'j_0}) < V(b_i, a_{ij_0})$, then we set $\bar{a}_{ij_0} = 0$.

(III) If $j_0 \in J_i$, $b_i \neq 0$ and there exists some $i' \in I$ ($i' \neq i$) such that $b_{i'} = 0$ and $a_{i'j_0} > 0$, then we set $\bar{a}_{ij_0} = 0$.

Then, two systems $S = \{x \in [0,1]^n : A\varphi x = b\}$ and $\bar{S} = \{x \in [0,1]^n : \bar{A}\varphi x = b\}$ have the same solutions.

**Definition 4.** For each $i \in I$, define $\bar{J}_i = \{j \in J : \bar{a}_{ij} \geq b_i\}$. Also, let $\bar{E}$ be the set of all vectors $e: I \to \cup_{i \in I} \bar{J}_i$ so that $e(i) = j \in \bar{J}_i, \forall i \in I$.

**Corollary 2.** $S \subseteq \{\underline{X}(e) : e \in \bar{E}\} \subseteq \{\underline{X}(e) : e \in E\}$.

**Proof.** Suppose that $\bar{a}_{ij_0} = 0$ by step (I). So, we have $j_0 \notin J_i$, i.e., $a_{ij_0} < b_i$. Thus, $\bar{a}_{ij_0} = 0 < b_i$ that means $j_0 \notin \bar{J}_i$. On the other hand, suppose that $\bar{a}_{ij_0} = 0$ by steps (II) or (III). So, $j_0 \in J_i$ ($a_{ij_0} \geq b_i$) and $b_i \neq 0$. In these cases, we have $\bar{a}_{ij_0} = 0 < b_i$, that is, $j_0 \notin \bar{J}_i$. Therefore, $\bar{J}_i \subseteq J_i$ that implies $\bar{E} \subseteq E$. Hence, $\{\underline{X}(e) : e \in \bar{E}\} \subseteq \{\underline{X}(e) : e \in E\}$. Now, the proof is resulted from relation (9) and Theorem 1. □

**Corollary 3.** $S = \cup_{e \in \bar{E}} [\underline{X}(e), \bar{X}]$.
**Proof.** The proof is resulted from relation (7) and Corollary 2. □

The following example illustrates the above-mentioned definitions.

**Example 1.** Consider the problem below with Dombi t-norm

$$\begin{bmatrix} 0.9452 & 0.4012 & 0.8976 & 0.6221 & 0.4368 & 0.8126 \\ 0.5271 & 0.1113 & 0.2456 & 0.3419 & 0.5271 & 0.2192 \\ 0.2073 & 0.8172 & 0.4386 & 0.4599 & 0.6152 & 0.2188 \\ 0.9111 & 0.7243 & 0.3274 & 0.8327 & 0.8327 & 0.5845 \end{bmatrix} \varphi \begin{bmatrix} 0.7243 \\ 0.5271 \\ 0.6152 \\ 0.8327 \end{bmatrix}$$

where

$$\varphi(x,y) = \begin{cases} 0 & x = 0 \text{ or } y = 0 \\ \left(1 + \sqrt{\left(\frac{1-x}{x}\right)^2 + \left(\frac{1-y}{y}\right)^2}\right)^{-1} & \text{otherwise} \end{cases}$$

(i.e., $\lambda = 2$). In this example, we have $J_1 = \{1,3,6\}$, $J_2 = \{1,5\}$, $J_3 = \{2,5\}$ and $J_4 = \{1,4,5\}$. According to relation (6), we attain $\hat{x}_1 = [0.7266 \ 1 \ 0.7335 \ 1 \ 1 \ 0.7675]$, $\hat{x}_2 = [1 \ 1 \ 1 \ 1 \ 1 \ 1]$, $\hat{x}_3 = [1 \ 0.6312 \ 1 \ 1 \ 1 \ 1]$ and $\hat{x}_4 = [0.8506 \ 1 \ 1 \ 1 \ 1 \ 1]$.
Therefore,
$\bar{X} = \min_{i \in I}\{\hat{x}_i\} = [0.7266 \ 0.6312 \ 0.7336 \ 1 \ 1 \ 0.7675]$
Also, based on Definition 1, we have

$\check{x}_1(1) = [0.7266 \ 0 \ 0 \ 0 \ 0 \ 0]$, $\check{x}_1(3) = [0 \ 0 \ 0.7335 \ 0 \ 0 \ 0]$, $\check{x}_1(6) = [0 \ 0 \ 0 \ 0 \ 0 \ 0.7675]$
$\check{x}_2(1) = [1 \ 0 \ 0 \ 0 \ 0 \ 0]$, $\check{x}_2(5) = [0 \ 0 \ 0 \ 0 \ 1 \ 0]$
$\check{x}_3(2) = [0 \ 0.6312 \ 0 \ 0 \ 0 \ 0]$, $\check{x}_3(5) = [0 \ 0 \ 0 \ 0 \ 1 \ 0]$
$\check{x}_4(1) = [0.8506 \ 0 \ 0 \ 0 \ 0 \ 0]$, $\check{x}_4(4) = [0 \ 0 \ 0 \ 1 \ 0 \ 0]$, $\check{x}_4(5) = [0 \ 0 \ 0 \ 0 \ 1 \ 0]$

The cardinality of set $E$ is equal to $|E| = \prod_{i \in I} |J_i| = 3 \times 2 \times 2 \times 3 = 36$. So, we have 36 solutions $\underline{X}(e)$ associated to 36 vectors $e$. For example, for $e = [1,5,5,5]$, we obtain $\underline{X}(e) = \max\{\check{x}_1(1), \check{x}_2(5), \check{x}_3(5), \check{x}_4(5)\}$ from Definition 1, 2 and 3 that means $\underline{X}(e) = [0.7266 \ 0 \ 0 \ 0 \ 1 \ 0]$. Now, from step (I) of Theorem 1, we have

$\bar{a}_{12} = \bar{a}_{13} = \bar{a}_{14} = \bar{a}_{15} = \bar{a}_{16} = \bar{a}_{22} = \bar{a}_{23} = \bar{a}_{24} =$

$\bar{a}_{26} = \bar{a}_{31} = \bar{a}_{33} = \bar{a}_{34} = \bar{a}_{36} = \bar{a}_{42} = \bar{a}_{43} = \bar{a}_{46} = 0$

In all of these cases, $a_{ij} < b_i$, that is, $j \notin J_i$. Also, from part (II), we can set $\bar{a}_{21} = \bar{a}_{41} = 0$. For example, $a_{21} = b_2$ (i.e., $1 \in J_2$), $b_2 \neq 0$, $a_{11} > b_1$ (i.e., $1 \in J_1$), $b_1 \neq 0$ and $0.7266 = V(b_1, a_{11}) < V(b_2, a_{21}) = 1$.

Hence, we have

$$\bar{A} = \begin{bmatrix} 0.9452 & 0.0000 & 0.8976 & 0.0000 & 0.0000 & 0.8126 \\ 0.0000 & 0.0000 & 0.0000 & 0.0000 & 0.5271 & 0.0000 \\ 0.0000 & 0.8172 & 0.0000 & 0.0000 & 0.6152 & 0.0000 \\ 0.0000 & 0.0000 & 0.0000 & 0.8327 & 0.8327 & 0.0000 \end{bmatrix}$$

and therefore, $\bar{J}_1 = \{1,3,6\}$, $\bar{J}_2 = \{5\}$, $\bar{J}_3 = \{2,5\}$ and $\bar{J}_4 = \{4,5\}$ that imply $|\bar{E}| = \prod_{i \in I} |\bar{J}_i| = 3 \times 1 \times 2 \times 2 = 12$.

## 3. LATTICIZED PROGRAMMING PROBLEM AND OPTIMAL SOLUTION

**Theorem 2.** If $S \neq \emptyset$, there exists a minimal solution of $S$, i.e. $\underline{x}^* \in \underline{S}$, such that $\underline{x}^*$ is an optimal solution of problem (1).

**Proof.** Let $z(x) = \max\{x_1, x_2, \ldots, x_n\}$. Furthermore, suppose that $z(\underline{x}^*) = \min\{z(\underline{x}) : \underline{x} \in \underline{S}\}$ where $\underline{x}^*$ is a minimal solution. Based on Corollary 3, for each $x' \in S$ there exist some $\underline{x} \in \underline{S}$ such that $\underline{x} \leq x'$, i.e., $\underline{x}_j \leq x'_j, \forall j \in J$. So, we have $\max\{\underline{x}_1, \underline{x}_2, \ldots, \underline{x}_n\} \leq \max\{x'_1, x'_2, \ldots, x'_n\}$ that implies $z(\underline{x}) \leq z(x')$. But, $z(\underline{x}^*) \leq z(\underline{x})$ which implies $z(\underline{x}^*) \leq z(x'), \forall x' \in S$. □

By combination of Theorem 2 and Corollary 3, it turns out that the optimal solution of problem (1) must be a vector $\underline{X}(e^*)$ for some $e^* \in \bar{E}$. Based on this fact, we can find the optimal solution of problem (1) by pairwise comparison between the elements of set $\{\underline{X}(e): e \in \bar{E}\}$. We now summarize the preceding discussion as an algorithm.

**Algorithm 1**
Given problem (1):
1. If $\bar{X} \notin S$, then $S$ is empty (Corollary 1).
2. Find solutions $\underline{X}(e), \forall e \in \bar{E}$ (Definition 4).
3. Find the minimal solutions, $\underline{S}$ by the pairwise comparison between the solutions $\underline{X}(e)$ (Corollary 3).
4. Find the optimal solution $\underline{X}(e^*)$ for problem (1) by the pairwise comparison between the objective values of the elements of $\underline{S}$ (Theorem 2).

## 3. NUMERICAL EXAMPLE

Based on the theories we built in previous sections, here we propose an algorithm for finding an optimal solution of problem (1). Consider the problem expressed in example 1 with the objective function $z(x) = \max\{x_1, x_2, \ldots, x_6\}$. According to Example 1, we calculated the maximum solution as

$\bar{X} = [0.7266 \quad 0.6312 \quad 0.7336 \quad 1 \quad 1 \quad 0.7675]$

Since $\bar{X} \in S$, then Corollary 1 implies that $S$ is nonempty.
The 12 vectors $\underline{X}(e)$ for each $e \in \bar{E}$ (recall that $|\bar{E}| = 12$) are calculated as follows:

$e_1 = [1\ 5\ 2\ 4]$
$\underline{X}(e_1) = [0.7266\ 0.6312\ 0\ 1\ 1\ 0]$

$e_2 = [3\ 5\ 2\ 4]$
$\underline{X}(e_2) = [0\ 0.6312\ 0.7336\ 1\ 1\ 0]$

$e_3 = [6\ 5\ 2\ 4]$
$\underline{X}(e_3) = [0\ 0.6312\ 0\ 1\ 1\ 0.7675]$

$e_4 = [1\ 5\ 5\ 4]$
$\underline{X}(e_4) = [0.7266\ 0\ 0\ 1\ 1\ 0]$

$e_5 = [3\ 5\ 5\ 4]$
$\underline{X}(e_5) = [0\ 0\ 0.7336\ 1\ 1\ 0]$

$e_6 = [6\ 5\ 5\ 4]$
$\underline{X}(e_6) = [0\ 0\ 0\ 1\ 1\ 0.7675]$

$e_7 = [1\ 5\ 2\ 5]$
$\underline{X}(e_7) = [0.7266\ 0.6312\ 0\ 0\ 1\ 0]$

$e_8 = [3\ 5\ 2\ 5]$
$\underline{X}(e_8) = [0\ 0.6312\ 0.7336\ 0\ 1\ 0]$

$e_9 = [6\ 5\ 2\ 5]$
$\underline{X}(e_9) = [0\ 0.6312\ 0\ 0\ 1\ 0.7675]$

$e_{10} = [1\ 5\ 5\ 5]$
$\underline{X}(e_{10}) = [0.7266\ 0\ 0\ 0\ 1\ 0]$

$e_{11} = [3\ 5\ 5\ 5]$
$\underline{X}(e_{11}) = [0\ 0\ 0.7336\ 0\ 1\ 0]$

$e_{12} = [6\ 5\ 5\ 5]$
$\underline{X}(e_{12}) = [0\ 0\ 0\ 0\ 1\ 0.7675]$

However, by the pairwise comparison, it is found that the three following vectors are the only minimal solutions of the feasible region:

$\underline{X}(e_{10}) = [0.7266\ 0\ 0\ 0\ 1\ 0]$
$\underline{X}(e_{11}) = [0\ 0\ 0.7336\ 0\ 1\ 0]$
$\underline{X}(e_{12}) = [0\ 0\ 0\ 0\ 1\ 0.7675]$

Since $z(\underline{X}(e_{10})) = z(\underline{X}(e_{11})) = z(\underline{X}(e_{12}))$, we conclude that the problem has three optimal points $\underline{X}(e_{10})$, $\underline{X}(e_{11})$ and $\underline{X}(e_{12})$ with the optimal value equal to one.

4. CONCLUSIONS

Considering the practical applications of the max-Dombi fuzzy relational equations in FRE theory and that of the latticized programming, a nonlinear optimization problem was studied with the maximum components function as the objective function subjected to the Dombi FRE. Since a system of the Dombi FRE is a non-convex set, an algorithm was presented to find an optimal solution by using the structural properties of the problem. For this purpose, a necessary and sufficient feasibility condition was firstly derived and then, the feasible region was completely determined in terms of one maximum and a finite number of minimal

دومین کنفرانس بین المللی یافته های پژوهشی مهندسی برق، کامپیوتر و مکانیک

2nd International Conference in Electrical, Computer and Mechanical Engineeringsolutions. It is proved that we can find the exact optimal solution of the proposed problem from the minimal solutions of the constraints, i.e., a system of max-Dombi FRE. Additionally, a numerical example was given to illustrate the presented algorithm.

5. **REFERENCES**

1. Chaira, Tamalika. "Fuzzy mathematical morphology using triangular operators and its application to images." Journal of Intelligent & Fuzzy Systems 28,5 (2015): 2269-2277.

2. Chang, Cheung-Wen, and Bih-Sheue Shieh. "Linear optimization problem constrained by fuzzy max–min relation equations." Information Sciences 234 (2013): 71-79.

3. Chen, Jiqian, and Jun Ye. "Some single-valued neutrosophic Dombi weighted aggregation operators for multiple attribute decision-making." Symmetry 9,6 (2017): 82.

4. Dombi, Jozsef. "On a certain class of aggregative operators." Information Sciences 245 (2013): 313-328.

5. D. Dubois, H. Prade, Fundamentals of Fuzzy Sets, Kluwer, Boston, 2000.

6. Fan, Y. R., Gordon H. Huang, and A. L. Yang. "Generalized fuzzy linear programming for decision making under uncertainty: Feasibility of fuzzy solutions and solving approach." Information Sciences 241 (2013): 12-27.

7. Freson, Steven, Bernard De Baets, and Hans De Meyer. "Linear optimization with bipolar max–min constraints." Information Sciences 234 (2013): 3-15.

8. Ghodousian, Amin. "Optimization of linear problems subjected to the intersection of two fuzzy relational inequalities defined by Dubois-Prade family of t-norms." Information Sciences 503 (2019): 291-306.

9. Ghodousian, Amin, and Fatemeh Elyasimohammadi. "On the optimization of Dombi non-linear programming." Journal of Algorithms and Computation 52,1 (2020): 1-36.

10. Ghodousian, Amin, and Maryam Raeisian Parvari. "A modified PSO algorithm for linear optimization problem subject to the generalized fuzzy relational inequalities with fuzzy constraints (FRI-FC)." Information Sciences 418 (2017): 317-345.

11. Ghodousian, Amin, Marjan Naeeimi, and Ali Babalhavaeji. "Nonlinear optimization problem subjected to fuzzy relational equations defined by Dubois-Prade family of t-norms." Computers & Industrial Engineering 119 (2018): 167-180.

12. Ghodousian, Amin, and Ali Babalhavaeji. "An efficient genetic algorithm for solving nonlinear optimization problems defined with fuzzy relational equations and max-Lukasiewicz composition." Applied Soft Computing 69 (2018): 475-492.

13. Li, Pingke, and Yuhan Liu. "Linear optimization with bipolar fuzzy relational equation constraints using the Łukasiewicz triangular norm." Soft Computing 18,7 (2014): 1399-1404.

14. Lin, Jun-Lin, Yan-Kuen Wu, and Sy-Ming Guu. "On fuzzy relational equations and the covering problem." Information Sciences 181,14 (2011): 2951-2963.

15. Liu, Chia-Cheng, Yung-Yih Lur, and Yan-Kuen Wu. "Linear optimization of bipolar fuzzy relational equations with max-Łukasiewicz composition." Information Sciences 360 (2016): 149-162.

16. Pedrycz W. Granular Computing: Analysis and Design of Intelligent Systems, CRC Press, Boca Raton. (2013).

17. Sanchez, Elie. "Solutions in composite fuzzy relation equations: application to medical diagnosis in Brouwerian logic." in M.M. Gupta. G.N. Saridis, B.R. Games (Eds.), Fuzzy Automata and Decision Processes, North-Holland, New York,(1977): pp. 221-234.
7